\documentclass[12pt]{amsart} \usepackage{amscd}

\newtheorem{theorem}{Theorem}%[section] (If you want theorem numbered
\theoremstyle{definition}

\theoremstyle{remark} 

\usepackage{amsfonts}% to get the /mathbb alphabet
\newcommand{\field}[1]{\mathbb{#1}}          \newcommand{\Q}{\field{Q}}
                   \newcommand{\Z}{\field{Z}}
                   \newcommand{\A}{\field{A}}

     \newcommand{\fp}{\mathfrak    p}
     
\newcommand{\fq}{\mathfrak q} \newcommand{\fa}{\mathfrak a}

\begin{document}
%\input{dspace11}
%\middlespace
%\end{document}

\title{On the gcd  of an infinite number of integers}

\author{T. N. Venkataramana}

\date{August 2, 2002}

\maketitle
\section*{Introduction}
\quad
%\end{document}

In  this  paper,  we  consider  the greatest  common  divisor  (to  be
abbreviated  gcd in  the  sequel) of  suitable  infinite sequences  of
integers  and  prove that  as  the  sequences  vary, the  gcd  remains
bounded.  These questions are motivated  by results of \cite{S} on the
congruence  subgroup problem.   Using  these results,  we also  obtain
bounds on  the indices of  certain congruence subgroups  of arithmetic
groups in algebraic groups. \\

The simplest example of the result on gcd is the following. Let $a$ and
$b$  be coprime positive  integers and  for an  integer $n\neq  0$ let
$\phi  (n)$ be  the number  of positive  integers not  exceeding $\mid
n\mid$  (the absolute  value of  $n$) and  coprime to  $\mid  n \mid$.
Consider the infinite sequence
$$\phi  (ax+b)~;  x=\cdots,-2,-1,0,1,2,\cdots  \leqno  (1)$$  and  let
$g(a,b)$  denote the  gcd of  the  numbers occurring  in the  sequence
$(1)$.  Then (see  Theorem 1 below), $g(a,b)$ is bounded  by $4$ , for
all  $a$  and $b$.   This  assertion  will be  shown  to  be a  simple
application of Dirichlet's  Theorem on the infinitude of  primes in an
arithmetic progression.\\

We  will also  prove an  analogue of  the above  result  for arbitrary
number fields.   Let $S$ be a finite  set of places of  a number field
$F$  containing all  the  archimedean  ones, $a$  and  $b$ be  coprime
elements in the  ring of $S$-integers $O_S$ , and  for a nonzero ideal
$\mathfrak  b$ of $O_S$,  $\phi({\mathfrak b})$  denote the  number of
units in  the residue class  ring $O_S/\mathfrak b$.  Write  $(y)$ for
the principal  ideal in $O_S$ generated  by the element  $y$ in $O_S$.
Then  we show  (see Theorem  2  below) that  the gcd  $g(a,b)$ of  the
infinite sequence of numbers
$$\phi  ((ax+b)); x\in  O_S, \leqno  (2)$$  is bounded  by a  constant
depending only on the number field  $F$.  In fact the constant will be
shown to  depend only on  the number of  roots of unity in  the number
field.   As in  the case  of $\Q$,  Theoerm 2  will be  shown to  be a
consequence of the Cebotarev density Theorem. \\

This result has an application to algebraic groups.  Theorem 2 will be
shown  to imply  that if  $\Gamma  $ is  an $S$-arithmetic  congruence
subgroup of a linear algebraic group $G$ defined over the number field
$F$, then the  congruence subgroups $\Gamma (ax+b)$ (of  $\Gamma $ for
$x$ varying in  $O_S$, generate a subgroup $\Gamma  '$ of finite index
$g$ in $\Gamma $, and the index $g$ is bounded by a constant depending
only on the  algebraic group $G$ and the number field  $F$ (and not on
the congruence subgroup $\Gamma $).\\

With these notations, we first prove (in $\S$1) the following:

\begin{theorem} Suppose that $a$ and $b$ are coprime rational integers
, and  let $\phi  $ denote  the Euler $\phi$  function.  Then  the gcd
$g=g(a,b)$ of the  sequence $\{ \phi(ax+b);~x\in \Z \}$  is bounded by
$4$; more precisely, $g$ has the following values.

\begin{enumerate}
\item[(1)]$g=1$ if $b\equiv \pm 1$ (mod~$a$) or $\pm 2$ (mod~$a$). 
\item[(2)]$g=4$ if $a$ is divisible by  four, $b$ is not as in Case 1,
$b=\equiv~1$ (mod~4) and $b$ is not a square (mod~$a$). 
\item[(3)]$g=2$ in all other cases.
\end{enumerate}
\end{theorem}

Theorem 1 is proved by  using Dirichlet's Theorem on the infinitude of
primes in arithmetic progressions. \\

Theorem 1 can be generalised to  arbitrary number fields. Let $F$ be a
number field and  $S$ a finite set of places of  $F$ including all the
archimedean ones.  Let $O_S$ denote the ring of $S$-$integers$ of $F$.
Given   a  non-zero   ideal   $\mathfrak{a}$  of   $O_S$,  let   $\phi
(\mathfrak{a})$ denote the number of  {\bf units} in the quotient ring
$O_S/\mathfrak{a}$ (of  course, the quotient has only  a finite number
of elements).   Given a non-zero  element $x\in O_S$, denote  by $\phi
(x)$ the  number $\phi ((x))$, where  $(x)$ is the  principal ideal of
$O_S$ generated by $x$. Recall that  two elements $a$ and $b$ of $O_S$
are said to be  coprime, if the ideal generated by $a$  and $b$ is the
whole ring $O_S$.  Let $\mu _F$  denote the group of roots of unity in
the number field $F$. \\

\begin{theorem} Let $a$ and $b$ denote two elements of $O_S$ which are
coprime. Let  $g$ be the  gcd of the  numbers $\phi (ax+b)$  where $x$
runs  through  elements of  $O_S$.   Then,  there  is a  constant  $C$
depending only  on the number field  $K$ (and not on  the elements $a$
and $b$) such that  $g\leq C$. In fact, $C$ may be  chosen to be $Card
(\mu _F)^2$. \\

More  precisely, $g$ divides  $Card (\mu  _F)^2$. That  is $v_l(g)\leq
v_l(Card  (\mu _F)  ^2)$ where  for a  prime $l$  and an  integer $m$,
$v_l(m)$ is the largest power of $l$ which divides $m$.\\

In particular,  $g$ is divisible only  by the primes  which divide the
order of the group $\mu _F$ of roots of unity in the number field $F$.
\end{theorem}

Theorem  2  has  an  application  to algebraic  groups.   Indeed,  its
formulation was suggested by a  problem in the arithmetic of algebraic
groups.  Let  $G$ be  a linear algebraic  group defined over  a number
field  $F$  and  let  $O_S$  be  as  before.   Let  $\Gamma  $  be  an
$S$-congruence subgroup of the group $G(O_S)$ of $O_S$-rational points
of $G$. Let  $a$ and $b$ be coprime elements in  $O_S$ and given $x\in
O_S$, let  $G(ax+b)$ be the  congruence subgroup corresponding  to the
principal ideal ($ax+b$) of $O_S$.   Let $\Gamma '$ be the subgroup of
$\Gamma  $ generated  by the  subgroups $\Gamma  \cap G(ax+b)$  as $x$
varies through elements of the ring $O_S$.\\

\begin{theorem} There exists an integer $g$ depending only on the
group  $G$ and  the number  field $F$  such that  for all  $\Gamma$ as
above, the index of $\Gamma '$ in $\Gamma $ is bounded by $g$.
\end{theorem}

\section{Proof of Theorem 1}

As we said before, the  proof is a repeated application of Dirichlet's
Theorem  on the infinitude  of primes  in arithmetic  progressions. We
will prove it by a (small) number of case by case checks. \\

[1]. Suppose that $b$ is congruent to either $\pm 1$ or $\pm 2$ modulo
$a$. Then, there exists an integer $x$ such that $ax+b$ is either $\pm
1$ or  $\pm 2$.   In either  case $\phi (ax+b)=1$;  since the  gcd $g$
divides each $\phi (ax+b)$, we obtain that $g=1$. \\

Let $l$ be a prime and  $m$ an integer. Define, as before, $v_l(m)$ to
be the largest integer $k$ such that $l^k$ divides $m$. \\

[2].  We first prove that
\[v_l(g)=0 \] 
if $l$ is an odd prime. Write $a'$ for the l.c.m. of $a$ and $l$.

If  $b$ is  coprime to  $l$ and  $b \not\equiv  1$ (mod~$l$),  then, by
Dirichlet's Theorem, the residue  class $b$ modulo $a'$ is represented
by  a  prime  $p$.   Thus,  there  exist integers  $x,  y$  such  that
$p=a'y+b=ax+b$, and $\phi  (ax+b)=\phi (p)=p-1\not\equiv 0$ (mod~$l$).
Therefore,  $g$  (which divides  $\phi  (p)$)  is  coprime to  $l$  as
well. That is $v_l(g)=0$.

If $b$ is not  coprime to $l$, then, by replacing $b$  by $a+b$ we may
still assume that  $b$ is coprime to $l$, as $a$  and $b$ are coprime.
Suppose now that $b\equiv  1$~(mod~$l$).  Since the group $(\Z/l\Z)^*$
of  units  in  the  ring  $\Z/l\Z$  has at  least  two  elements,  and
$(\Z/a'\Z)^*\rightarrow (\Z/l\Z)^*$ is surjective, we may write $b=uv$
modulo $a'$ with $u\not\equiv 1$  and $v\not\equiv 1$ modulo the prime
$l$. Moreover,  by Dirichlet's  Theorem, $u$ and  $v$ modulo  $a'$ are
represented by  (distinct) primes $p$  and $q$. Thus, there  exist $x$
and $y$  with $uv=a'y+b=ax+b$ and  $p\equiv u$ and $q\equiv  v$ modulo
$a'$.    Moreover,  $p-1\not\equiv  0$~(mod~$l$)   and  $q-1\not\equiv
0$~(mod  ~$l$).   Thus,  $\phi (pq)=(p-1)(q-1)\not\equiv  0$~(mod~$l$)
whence $v_l(g)=0$. This completes [2]. \\

[3]. Now we consider the case when $l$ is even ($l=2$). We show that
\[v_2(g)\leq 2.\]
{\bf (3.1)}. Suppose that $a$ is even and that $a^*$ is the lcm of $a$
and  $4$.  Now, $b$  is odd  since $b$  is coprime  to $a$.   Thus $b$
defines a class in $(\Z/4\Z)^*$.

If $b$ is  a non-trivial class in $(\Z/4\Z)^*$,  then, since the class
of $b$  modulo $a^*$ is represented  by a prime $p$,  one has $p=ax+b$
for  some $x$  with $\phi  (ax+b)=\phi (p)=p-1\neq  0$~(mod~$4$).  Thus,
$v_2(\phi (ax+b))\leq 1$ whence $v_2(g)\leq 1$.

If  $b$ {\it is}  the trivial  class modulo  $4$, then  write $b\equiv
(-b)(-1)$~(mod~$a^*$). Both $(-b)$ and  $(-1)$ are represented by primes
$p$ and  $q$ say. Then  $p\equiv q\equiv (-1)$~(mod~$4$),  whence, $p-1$
and $q-1$ are divisible at most by $2$. Thus, there exists an $x$ with
$ax+b=pq$,  with  $\phi  (ax+b)=\phi (p)\phi  (q)=(p-1)(q-1)\not\equiv
0$~(mod ~$8$). Hence $v_2(g)\leq v_2(\phi (ax+b))\leq 2$.\\

\noindent{\bf (3.2)}. 
Suppose $a$ is  odd. Let $a^*=4a$ denote the lcm of  $a$ and $4$.  The
group  $(\Z/a^*\Z)^*$  of  units  is  a product  of  $(\Z/a\Z)^*$  and
$(\Z/4\Z)^*$.  Hence there exists a  prime $p$ such that it represents
$b$ modulo  $a$ and $-1$  modulo $4$.  Thus  there exists an  $x$ with
$p=ax+b$.   Moreover $\phi  (ax+b)=\phi  (p)=p-1\not\equiv 0$~(mod~$4$),
whence $v_2(g)\leq 1$. \\

[4]. Suppose  that $ax+b$ is at  least $3$ in  absolute value.  Assume
also that $a\equiv 0 $ (mod~$4$), $b=1$ (mod~$4$) and that $b$ is not a
square  (mod $a$).  Write  $ax+b=p_1^{e_1}\cdots p_h^{e_h}$,  with $p_i$
distinct  primes and  $e_i$  positive integers.   Clearly, the  primes
$p_i$ are all odd.

If $h\geq  2$, then  since $\phi (ax+b)$  is divisible by  the product
$(p_i-1)\cdots (p_h-1)$, it follows that $\phi (ax+b)$ is divisible by
$2^h$ and hence by $4$.
 
If $h=1$  and $p_1=3~(mod~4)$, then  by passing to  congruences modulo
$4$ (note  that $4$ divides $a$  by assumption), we  see that $1\equiv
b\equiv p_1^{e_1}\equiv (-1)^{e_1}$ modulo  $4$, whence $e_1$ is even,
contradicting the assumption that $b$ is not a square modulo $a$.

If    $h=1$    and   $p_1\equiv    1$    (mod    $4$),   then    $\phi
(ax+b)=p_1^{e_1-1}(p_1-1)$ is divisible modulo $4$.

We have thus  proved part (2) of Theorem 1, that  for all $x$, $\phi
(ax+b)$ is divisible  by $4$; that is, $4$ divides  $g$. Since we have
proved in all cases that $v_2(g)\leq  2$, it follows that $g=4$ if the
conditions of Part (2) hold. \\

[5]. Turn now to the other parts of Theorem 1. If $a\equiv 2~(mod~4)$,
then write  $a=2a_1$, with $a_1$  odd.  By Dirichlet's  Theorem, there
exist  infinitely many  primes $p$  with $p\equiv  b$ (mod  $a_1$) and
$p\equiv 3$ (mod  $4$). Note that $b\equiv p\equiv  1$ (mod $2$) since
both  $b$ and  $p$ are  odd.  Thus, $p\equiv  b$ (mod  $a$) and  $\phi
(p)=p-1\not\equiv  0$ (mod $4$)  whence $v_2(g)\leq  1$. Since  $g$ is
always even in this case, it follows that $g=2$.\\

If $b$ is $\pm  1$ or $\pm 2$, then, there exists  an integer $x$ such
that  $ax+b=\epsilon$  with  $\epsilon  =\pm  1$ or  $\epsilon  =  \pm
2$. Clearly, then $\phi (\epsilon )= 1$ whence $g=1$. \\

If $a$ and  $b$ are as in Part  (3) of Theorem 1, then  $ax+b$ is at
least  3 in  absolute value  for all  $x$ and  hence  $\phi (ax+b)\geq
2$. If $a$ is divisible by $4$,  but $b$ is a square say $c^2$, modulo
$a$,  then $c$  is represented  by a  prime $p$  modulo $4$  which (by
replacing $c$ by $-c$ if necessary)  may be assumed to be congruent to
$-1$ modulo  $4$.  Hence  $p^2=ax+b$ for some  integer $x$,  and $\phi
(ax+b)=p(p-1)\not\equiv 0$(mod $4$).  Therefore, $v_2(g)\leq 1$ whence
$g=2$.  The rest of Theorem 1 may be proved in exactly the same way.\\

\section{Proof of Theorem 2.} 

\subsection{Preliminaries.} Denote by $\mu $ the order of the group of 
roots of unity in the number field  $F$.  Let $l$ be a prime and $l^e$
the largest  power of $l$  which divides $\mu  $.  Let $\omega $  be a
primitive  $l^{e+1}$-th  root of  unity  in  an  algebraic closure  of
$F$. The Galois group of the extension $F(\omega )/F$ is a subgroup of
the group $(\Z/l^{e+1}\Z)^*$ of  units and is non-trivial, since $l^e$
is the highest power of $l$  which divides $\mu $.  Moreover, if $\fp$
is a prime ideal in $O_S$ which is unramified in $F(\omega)$ such that
its Frobenius is  a non-trivial element in $Gal  (F(\omega )/F)$ then,
$N(\fp)-1$ is divisible by $l^e$ but not by $l^{e+1}$.\\

\subsection{Notation.} 
Denote by $F_a$ the ray class  field of $F$ corresponding to the ideal
$O_S$.   Thus, by definition,  $F_a/F$ is  an abelian  extension whose
Galois group  is isomorphic -under  the Artin Reciprocity map-  to the
quotient  $G(F_a/F) =\A_F^*/F^*(\prod  U_v)$. Here  for  a commutative
ring $R$, we  denote by $R^*$ the  group of units of $R$  and $A_F$ is
the ring  of adeles  of $F$.  For  each finite  place $v$ of  $F$, let
$O_v$  be the  maximal compact  subring of  the completion  $F_v$.  If
$v\in S$,  let $U_v=K_v^*$  and if $v\notin  S$, let $U_v$  denote the
subgroup of elements  of $O_v^*$ congruent to $1$  modulo the ideal of
$O_v$ generated  by $a$.  Note  that $S$ contains all  the archimedean
primes and  that if $v\notin S$ and  if $v$ does not  divide $a$, then
$U_v=O_v^*$. \\

Given an element  $z=\prod (z_v)_v$ of $A_F^*$, define  its {\bf norm}
by $\prod  (\mid z\mid  _v)$ for all  $v\notin S$. In  particular, its
restriction to $\prod O_v^*$ (where $v$ runs through all finite places
of $F$) is trivial, and hence  the norm is a homomorphism on the group
of fractional  ideals of $O_S^*$. Thus,  if $\fp$ is a  prime ideal in
$O_S$ and $v$ the corresponding non-archimedean absolute value on $F$,
then the norm is $1/Card (O_S/\fp)$.\\

By the weak approximation theorem, it is clear that the group
\[G(F_a/F)=\A_F^*/F^*(\prod U_v)\] 
is isomorphic to the group
\[\A_F(S)^*/O_S(a)^*(\prod O_v^*)\]
where $O_S(a)^*$  is the group  of units in  the ring $O_S$  which are
congruent to the identity modulo  the principal ideal $a$ of $O_S$ and
the product  is over  all the  places $v$ of  $F$ not  in $S$  and not
dividing $a$. \\

\subsection{Proof of Theorem 2.} Let $E$ be the compositum of the ray 
class field $F_a$ and the cyclotomic extension $F(\omega )/F$ ($\omega
$  being a  primitive $L^{e+1}$-th  root of  unity). The  Galois group
$G(E/F)$ is  abelian and has surjective  homomorphisms onto $G(F_a/F)$
and onto  $G(F(\omega )/F)$.  The  group $G(E/F)$ is  generated-by the
Cebotarev Density  Theorem-by the Frobenius  elements corresponding to
prime ideals $\fp$ of the ring  $O_S$ (which may also be assumed to be
unramified in the extension $E/F$).  In particular, the image of $b\in
F$ ($F^*$  thought of as  a subgroup of  $\A _F(S)^*$ (= the  group of
$S$-ideles))  in $G(E/F)$  is represented  by (infinitely  many) prime
ideals $\fp$.   This means that  $\fp=(ax+b)O_S$ for some  $x\in O_S$.
In particular, $\phi (\fp)=\phi (ax+b)$. \\

Case 1: If moreover, there  exists an element $\sigma \in G(E/F)$ such
that its image in $G(F_a/F)$ is represented by the element $b$ and its
image  in   $G(F(\omega  )/F)$   is  non-trivial,  then   $\sigma$  is
represented by a prime ideal $\fp$, with the following properties.  By
the choice of  $\sigma$, one gets an $x\in  O_S$ such that $ax+b=\fp$.
Thus  $\phi (ax+b)=\phi  (\fp)$  and $\fp$  has  non-trivial image  in
$G(F(\omega )/F)$.  Hence  $\phi (\fp)$ is divisible by  $l^e$ but not
by $l^{e+1}$ (see (2.1)).\\

Case 2: If  Case 1 is not possible, then this  means that the preimage
of  $b\in  G(F_a/F)$  in  $G(E/F)$  maps to  the  trivial  element  in
$G(F(\omega )/F)$. Fix any $\sigma$ in the preimage, and let $\tau \in
G(E/F)$ map  non-trivially into  $G(F(\omega )/F)$ (recall  from (2.1)
that the latter group is non-trivial).  One may write $\sigma =(\sigma
\tau)(\tau)^{-1}$. Note that both  $\sigma \tau$ and $\tau ^{-1}$ have
non-trivial images  in $G(F(\omega  )/F)$. Represent both  elements by
prime  ideals $\fp$  and  $\fq$ respectively.   Then  $\phi (\fp)$  is
divisible  by  $l^e$  but   not  by  $l^{e+1}$  (similarly  for  $\phi
(\fq)$). However, $\fp\fq=ax+b$ for  some $x\in O_S$. Therefore, $\phi
(ax+b)$ is  divisible by $l^{2e}$  but by no  higher power of  $l$. In
particular, $v_l(g)\leq 2e$.\\

The combination of Cases 1 and 2 proves Theorem 2.

\section{A Question}

In sections  1 and  2 we considered  the gcd  of $\phi (ax+b)$  as $x$
varies. Here $ax+b$  is a {\bf linear} polynomial  of degree one, with
coprime  coefficients.   We  now  replace  $  ax+b$  by  a  polynomial
$P(x)=a_nx^n+\cdots+a_0$  where   $a_0,\cdots,a_n$  are  integers  (or
S-integers)  which are ``coprime''  i.e. generate  the unit  ideal. We
have  the  following  question,   arrived  at  in  conversations  with
M. S. Raghunathan.

\subsection{Question:} Does There exist a constant $C=C(n)$ depending 
only on the degree $n$ of the  polynomial $P$ such that the gcd of the
integers $\phi  (P(x)): x=0,1,\cdots $  is bounded by $C(n)$  for {\bf
all} polynomials $P$ of degree $n$ with coprime coefficients ?\\

In  case the number  field $F$  is different  from $\Q$,  the question
amounts to asking if the constant $C$ depends only on the number field
$F$ and on the finite set $S$ of places.

\section{Relation to the Congruence Subgroup Problem} 

Consider the group $SL_2$ over a number field $F$. Suppose the set $S$
is such that  the group $O_S^*$ of units is  infinite.  Given an ideal
$\fa$ of $O_S$, let $G(\fa )$ denote the principal congruence subgroup
of  $SL_2(O_S)$  of  level $\fa$:  this  is  the  set of  matrices  of
determinant one, with entries in the ring $O_S$ which are congruent to
the identity matrix  modulo the ideal $\fa$. Let  $E(\fa )$ denote the
normal  subgroup  of $SL_2(O_S)$  generated  by  the upper  triangular
matrices with  $1$'s on  the diagonal and  which are congruent  to $1$
modulo the ideal $\fa$. Then, $E(\fa )$ is a normal subgroup of $G(\fa
)$. \\

In the  course of  proof of the  congruence subgroup property  for the
group  $SL_2(O_S)$, Serre (in  \cite{S}) considers  the action  of the
group $T(O_S)\simeq O_S^*$ on  the quotient $G(\fa )/E(\fa)$. Here $T$
is the group of diagonals in  $SL_2$.  Serre shows that this action is
trivial on the  subgroup of $T(O_S)$ generated by  $\mu $-th powers of
elements of  $T(O_S)$ (recall that  $\mu $ is  the number of  roots of
unity in  the number  field $F$).  This  is achieved by  proving first
that if a  coset class $\xi $ in the quotient  $SL_2(\fa )/E(\fa )$ is
represented by a matrix $A=\left (\begin{array}{cc} a & b \\ c & d
\end{array} \right)$, then, the group $T(aO_S)$ acts trivially on the
element $\xi  $.  Next, if  the element $A$  is replaced by  $BA$ with
$B=\left (\begin{array}{cc} 1 & x \cr 0 & 1 \end{array} \right )$ with
$x\in \fa $, then $BA=\left( \begin{array}{cc} a' & b' \\ c' &d'
\end{array} \right)$ with $a'=ax+b$.  Consequently, even the group
$T((ax+b)O_S)$ acts  trivially on the  element $\xi $. Thus  the group
$T_0$ generated  by the subgroups $T((ax+b)O_S)$, as  $x$ varies, acts
trivially on the  element $\xi$. By arguments similar  to the proof of
Theorem 2  (indeed, the  proof of  Theorem 2 in  the present  paper is
modelled on  \cite{S}), Serre then shows  that the group  of $\mu $-th
powers  of $T(O_S)$ is  contained in  the group  $T_0$ and  hence acts
trivially on $\xi$.\\

We now state a generalisation of this result. Let $G\subset GL_n$ be a
linear algebraic group defined over a  number field $F$.  Let $S$ be a
finite set  of places of $F$  containing all the  archimedean ones and
let $O_S$  denote the  ring of $S$-integers  of $F$. Given  a non-zero
ideal  $\fa $  in $O_S$,  denote by  $G(\fa)$ the  subgroup  of $G\cap
GL_n(O_S)$  whose entries  are congruent  to the  identity  modulo the
ideal $\fa $. \\

Suppose now  that $a$ and $b$  are coprime elements in  $O_S$ and that
$x\in O_S$. Denote, as  before, by $G(ax+b)$, the principal congruence
subgroup of level $ax+b $.  Let $\Gamma \subset G(F)$ be an arithmetic
subgroup  and  let $\Gamma  _{a,b}$  be  the  group generated  by  the
subgroups $G(ax+b)\cap \Gamma$ for all $x\in O_S$ of $G(O_S)$.

\begin{theorem} \label{index} 
With the foregoing  notation, the group $\Gamma _{a,b}$  has its index
in $\Gamma $  bounded by a constant dependent  only on $G\subset GL_n$
and the number field $F$,  but neither on the arithmetic group $\Gamma
$ nor on the elements $a$ and $b$.
\end{theorem}

We begin  by noting that  for the multiplicative group  $G={\bf G}_m$,
this is exactly  Serre's Theorem: the index is bounded  by $\mu $, the
number of roots of unity in the number field.  \\

The  proof of Theorem  \ref{index} is  easy when  $G$ is  unipotent or
semi-simple, since one may appeal to the strong approximation theorem,
a version  of which  implies that  if $\fa $  and ${\mathfrak  b}$ are
coprime  ideals, then  $G(\fa )$  and $G({\mathfrak  b} )$  generate a
subgroup  of index  bounded by  a  constant independent  of $\fa$  and
${\mathfrak    b}$   (see    \cite{R},   section    2    for   similar
considerations). Thus, the proof of Theorem 3 may be reduced easily to
that for tori, and hence, by  the structure theorem for tori , to that
for the  groups ${\bf G}_m$ over  some {\it finite  extensions} $E$ of
$F$.  Thus,  we are  reduced to  the case of  ${\bf G}_m$  over number
fields, for which Serre's result is already available.

\section{Acknowledgements} The author is grateful to 
S.D.Adhikari,  R.Balasubramanian  and  to  K.Srinivas for  their  kind
invitation  to take part  and to  contribute to  the proceedings  of a
conference  in   honour  of  Prof.Subba  Rao  held   at  Institute  of
Mathematical Sciences, Chennai. He is also grateful to the referee for
a careful reading of the MS, for making several necessary corrections,
and for  pointing out a  gap (which is  now filled!)  in the  proof of
Theorem 1. \\

The author acknowledges  the help of M. S.  Raghunathan in arriving at
the question raised  in section 3. He is also  grateful to J.-P. Serre
for suggesting a possible solution through sieve theoretic methods.


\begin{thebibliography}{JPSH}

\bibitem [S]{S} J.-P.Serre, Le Probleme de groupes de congruence pour
  $SL_2$, Ann. Math. {\bf 92} (1970), 489-527.\\

\bibitem [R] {R} M.S.Raghunathan, On the Congruence Subgroup Problem II,
  Invent.Math., {\bf 85} (1986), 73-117. \\

\bibitem [V] {V} T.N.Venkataramana, On Systems of Generators of
  Higher Rank Arithmetic Subgroups, Pac. J. Math., {\bf 166}, No.1
  (1994), 193-212. 

\end{thebibliography}
\end{document}